\documentclass[12pt,a4paper,draft]{article}
\usepackage[english,russian]{babel}
\usepackage{amssymb, amsmath, latexsym, amsfonts, amsthm}
\begin{document}
\renewcommand{\refname}{References}
\newtheorem{theorem}{Theorem}
\newtheorem*{theorem A}{Theorem A}
\newtheorem{lemma}{Lemma}
\newtheorem{corollary}{Corollary}
\newtheorem*{proposition}{Proposition}
\begin{center}
On some classes of infinitely differentiable functions
\end{center}
\begin{center}
I.Kh. Musin
\end{center}

\vspace {0.1cm}

\renewcommand{\abstractname}{}
\begin{abstract}
{\sc Abstract}.
Spaces of infinitely differentiable functions on ${\mathbb R}^n$ (more general than Gelfand-Shilov spaces of type $W_M$) are considered in the article. For this space Paley-Wiener type theorem is obtained. 

\vspace {0.25cm}






\vspace {0.25cm}

\end{abstract}

\vspace {0.3cm}
 
\section{Introduction}

{\bf 1.1. On the aims of the article}. In \cite {CO} there were introduced and explored weighted spaces of entire functions which are natural generalizations of the Gelfand-Shilov $W^{\Omega}$-type spaces \cite {GS1}, \cite {GS2}.     
Under some conditions on weights the Fourier transforms of functions belonging to these spaces were described. Thus,  more general than Gelfand-Shilov $W_M$-type spaces appeared in \cite {CO} too. 
Here we propose a direct way of generalization of the Gelfand-Shilov $W_M$-type spaces and study the Fourier transforms of functions belonging to such spaces with aim to make some results from \cite {CO} more convenient for applications.

Let ${\mathcal H} = \{h_{\nu}\}_{{\nu}=1}^{\infty}$ be  a family of convex functions $h_{\nu}\colon {\mathbb R}^n \to {\mathbb R}$ 
such that for each $\nu \in {\mathbb N}$:

1) $h_{\nu}(x) = h_{\nu}(\vert x_1 \vert, \ldots , \vert x_n \vert), \ x = (x_1, \ldots , x_n)\in {\mathbb R}^n$;

2) the restriction of $h_{\nu}$ to $[0, \infty)^n$ is nondecreasing in each variable;

3) $\displaystyle \lim_{x \to \infty} \frac {h_{\nu}(x)}{\Vert x \Vert}= + \infty$;

4) for each $M>0$ there exists a constant $A_{\nu, M}>0$ such that
$$
h_{\nu}(x) \le 
\sum \limits_{1 \le j \le n: x_j \ne 0} x_j \ln\frac {x_j}{M} + A_{\nu, M}, \ x = (x_1, \ldots , x_n) \in [0, \infty)^n;
$$ 

5) $
h_{\nu}(x) - h_{\nu + 1}(x) \ge \ln 2 \cdot \sum \limits_{j=1}^n x_j - \gamma_{\nu}, \ 
x = (x_1, \ldots , x_n) \in  [0, \infty)^n;
$

6) $h_{\nu + 1}(x + y) \le h_{\nu}(x) + h_{\nu}(y) + l_{\nu}, \ x, y \in  [0, \infty)^n$.

For each $\nu \in {\mathbb N}$ and $m \in {\mathbb Z}_+$ define the normed space
$$
G_m(h_{\nu}) = \{f \in C^m({\mathbb R}^n): 
\Vert f \Vert_{m, h_{\nu}}
 = \sup_{x \in {\mathbb R}^n, \beta \in {\mathbb Z}_+^n, \atop \vert \alpha \vert \le m}  
\frac 
{\vert x^{\beta}(D^{\alpha}f)(x) \vert}
{\beta! e^{-h_{\nu}(\beta)}} < \infty \}.
$$
Let $G(h_{\nu})= \bigcap \limits_{m = 0}^{\infty} G_m(h_{\nu})$.
Endow $G(h_{\nu})$ with the topology defined by the family of norms 
$\Vert \cdot \Vert_{m, h_{\nu}}$ ($m \in {\mathbb Z}_+$).
Considering $G({\mathcal H}) = \bigcup \limits_{\nu = 1}^{\infty} G(h_{\nu})$ 
supply it with the topology of the inductive limit of spaces $G(h_{\nu})$.

For $f \in G({\mathcal H})$ let 
$$
\hat f(x) = \frac {1}{(\sqrt {2 \pi})^n} \int_{{\mathbb R}^n} f(\xi) e^{i \langle x, \xi \rangle} \ d \xi , \ x \in {\mathbb R}^n,
$$
be its Fourier transform.

One of the aims of the article is to describe the space of the Fourier transforms of functions belonging to $G({\mathcal H})$.

{\bf 1.2. Some notations}.  
For $u=(u_1, \ldots , u_n), v=(v_1, \ldots , v_n) \in {\mathbb R}^n \ ({\mathbb C}^n)$ let denote
$\langle u, v \rangle  = u_1 v_1 + \cdots + u_n v_n$ and let $\Vert u \Vert$  be the Euclidean norm of $u$. 

For $\alpha = (\alpha_1, \ldots , \alpha_n) \in {\mathbb Z}_+^n$, 
$x =(x_1, \ldots , x_n) \in {\mathbb R}^n$, 
$z =(z_1, \ldots , z_n) \in {\mathbb C}^n$ let
$\vert \alpha \vert = \alpha_1 + \ldots  + \alpha_n$, 
$\alpha! = \alpha_1! \cdots \alpha_n!$, 
$x^{\alpha} = x_1^{\alpha_1} \cdots x_n^{\alpha_n}$, $z^{\alpha} = z_1^{\alpha_1} \cdots z_n^{\alpha_n}$, 
$D^{\alpha}=
\frac {{\partial}^{\vert \alpha \vert}}{\partial x_1^{\alpha_1} \cdots \partial x_n^{\alpha_n}}$ .

For 
$\alpha = ({\alpha}_1, \ldots , {\alpha}_n)$ and 
$\beta = ({\beta}_1, \ldots , {\beta}_n) \in {\mathbb Z_+^n}$ the notation
$\alpha \le \beta $ indicates that
${\alpha}_j \le {\beta}_j$ ($j = 1, 2, \ldots , n$) 
and in such case $\binom {\beta}{\alpha}:= \prod \limits_{j=1}^{n} \binom {\beta_j}{\alpha_j}$, 
where $\binom {\beta_j}{\alpha_j}$
are the binomial coefficients. 

By $s_n$ denote the surface area of the unit sphere in ${\mathbb R}^n$.
 
${\mathbb R}_+^n:= [0, \infty)^n$.
For $a > 0$ let $\tilde a = (a,  \ldots , a) \in {\mathbb R}_+^n$.

If $[0, \infty)^n \subseteq X \subset {\mathbb R}^n$ then for a function $u$ on $X$ denote by $u[e]$ the function defined by the rule:
$u[e](x) = u(e^{x_1}, \ldots, e^{x_n}), \ x = (x_1, \ldots , x_n) \in  {\mathbb R}^n$.

The Young-Fenchel conjugate of a function $g:{\mathbb R}^n \to [-\infty, + \infty]$ is the function 
$g^*:{\mathbb R}^n \to [-\infty, + \infty]$ defined by
$$
g^*(x) = \displaystyle \sup \limits_{y \in {\mathbb R}^n}(\langle x, y \rangle - g(y)), \ x \in {\mathbb R}^n. 
$$

Denote the space of entire functions on ${\mathbb C}^n$ by $H({\mathbb C}^n)$.  

The Fourier transformation $f \in G({\mathcal H}) \to  \hat f$ is denoted by ${\mathcal F}$.

{\bf 1.3. Main results}. For $\nu \in {\mathbb N}$, $m \in {\mathbb Z}_+$ consider the normed space 
$$
{\mathcal E}_m(h_{\nu}) =\{f \in  C^{\infty}({\mathbb R}^n): 
\rho_{m, \nu}(f) = \sup_{x \in {\mathbb R}^n, \alpha \in {\mathbb Z}_+^n} 
\frac {(1+ \Vert x \Vert)^m \vert (D^{\alpha}f)(x) \vert}{\alpha! e^{-h_{\nu}(\alpha)}} < \infty \}.
$$
Let ${\mathcal E}(h_{\nu}) = \bigcap \limits_{m=0}^{\infty}{\mathcal E}_m(h_{\nu})$.
Endow ${\mathcal E}(h_{\nu})$ with the topology defined by the family of norms $\rho_{m, \nu}$ ($m \in {\mathbb Z}_+$). 
Let ${\mathcal E}({\mathcal H}) = \bigcup \limits_{\nu=1}^{\infty}{\mathcal E}(h_{\nu})$. 
Supply ${\mathcal E}({\mathcal H})$ with an inductive limit topology of spaces ${\mathcal E}(h_{\nu})$. 

Note that each function $f \in {\mathcal E}({\mathcal H})$ admits (the unique)  
extension to entire function in ${\mathbb C}^n$. 
Indeed, in view of the condition 3) on ${\mathcal H}$ for each $\varepsilon >0$ there is a constant $c_{\varepsilon}(f)>0$ such that  
for all 
$\alpha \in {\mathbb Z}_+^n$ we have that
\begin{equation}
\vert (D^{\alpha}f)(x)\vert \le c_{\varepsilon}(f) {\varepsilon}^{\vert \alpha \vert}\alpha!, \ x \in {\mathbb R}^n. 
\end{equation}
Hence, the sequence
$(\sum \limits_{\vert \alpha \vert \le k} \frac 
{(D^{\alpha}f)(0)}{\alpha!} x^{\alpha})_{k=1}^{\infty}$ 
converges to $f$ uniformly on compacts of ${\mathbb R}^n$.  
Also from (1) it follows that the series 
$
\displaystyle \sum_{\vert \alpha \vert \ge 0} \frac {(D^{\alpha}f)(0)}{\alpha!} z^{\alpha}
$
converges uniformly on compacts of ${\mathbb C}^{n}$. Hence, 
its sum $F_f(z)$ is an entire function in ${\mathbb C}^n$. 
Note that $F_f(x) = f(x), \ x \in {\mathbb R}^n$. 
Denote by ${\mathcal A}$ the mapping $f \in {\mathcal E}({\mathcal H}) \to F_f$. 

For each $\nu \in {\mathbb N}$ define a function $\varphi_{\nu}$ on ${\mathbb R}^n$ by
$$
\varphi_{\nu}(x) = h_{\nu}^*(\ln (1 + \vert x_1 \vert), \ldots , \ln (1 + \vert x_n \vert)), \ 
x = (x_1, \ldots , x_n) \in  {\mathbb R}^n.
$$
Note that in view of condition 4) on ${\mathcal H}$ for each $\nu \in {\mathbb N}$ we have that
$$
\displaystyle \lim_{x \to \infty} \frac {\varphi_{\nu}(x)}{\Vert x \Vert}= + \infty.
$$
Obviously, the restrictions of functions $\varphi_{\nu}$ to $[0, \infty)^n$ are nondecreasing in each variable.
 
Next, for each $\nu \in {\mathbb N}$ and $m \in {\mathbb Z}_+$ consider the normed space
$$
E_m(\varphi_{\nu}) = \{f \in H({\mathbb C}^n): p_{\nu, m}(f) = \sup_{z \in {\mathbb C}^n} 
\frac 
{\vert f(z)\vert (1 + \Vert z \Vert)^m}
{e^{\varphi_{\nu} (Im \, z)}} < \infty \}.
$$
Obviously, $E_{m+1}(\varphi_{\nu})$ is continuously embedded in $E_m(\varphi_{\nu})$. 
Let $E(\varphi_{\nu})= \bigcap \limits_{m=0}^{\infty} E_m(\varphi_{\nu})$.
Endow $E(\varphi_{\nu})$ with a projective limit topology of spaces $E_m(\varphi_{\nu})$. 
Note that if $f \in E(\varphi_{\nu})$ then  
for each $m \in {\mathbb Z}_+$ (using the condition 5)) 
$$
p_{\nu+1, m}(f) \le e^{\gamma_{\nu}} p_{\nu, m}(f).
$$ 
Hence, $E(\varphi_{\nu})$ is continuously embedded in $E(\varphi_{\nu + 1})$ for each $\nu \in {\mathbb N}$. 
Let $\varPhi =\{\varphi_{\nu}\}_{{\nu}=1}^{\infty}$ and $E(\varPhi)= \bigcup \limits_{\nu=1}^{\infty} E(\varphi_{\nu})$. With the usual operations of addition and multiplication by complex numbers 
$E(\varPhi)$
is a linear space. 
Supply $E(\varPhi)$ with a topology of the inductive limit of spaces $E(\varphi_{\nu})$.

The following results will be proved in the article.

\begin{theorem} 
The mapping ${\mathcal A}$ establishes an isomorphism between the spaces ${\mathcal E}({\mathcal H})$ and 
$E(\varPhi)$.
\end{theorem}

\begin{proposition}
${\mathcal F}$ is an injective linear continuous operator  
from $G({\mathcal H})$ into ${\mathcal E}({\mathcal H})$.
\end{proposition} 

\begin{theorem}  
The mapping ${\mathcal A} {\mathcal F}$ 
establishes an isomorphism between the spaces $G({\mathcal H})$ and 
$E(\varPhi)$.
\end{theorem}

These results are proved in sections 3, 4  and 5. The direct consequence of them is the following theorem.

\begin{theorem}
The Fourier transformation establishes an isomorphism between $G({\mathcal H})$ and 
${\mathcal E}({\mathcal H})$.
\end{theorem}

Further, 
let ${\mathcal M} = \{{M_{\nu}}\}_{\nu=1}^{\infty}$ be a family of convex functions 
$M_{\nu}: {\mathbb R}^n \to {\mathbb R}$ such that for each $\nu \in {\mathbb N}$:

$j_1)$ $M_{\nu}(x) = M_{\nu}(\vert x_1 \vert, \ldots , \vert x_n \vert), \ x = (x_1, \ldots , x_n)\in {\mathbb R}^n$;

$j_2)$ the restriction of $M_{\nu}$ to $[0, \infty)^n$ is nondecreasing in each variable;

$j_3)$ $\displaystyle \lim_{x \to \infty} \frac {M_{\nu}(x)}{\Vert x \Vert}= + \infty$.

$j_4)$ there exists a constant $A_{\nu}>0$ such that
$$
M_{\nu + 1}(2x) \le M_{\nu} (x) + A_{\nu}, \ x \in {\mathbb R}^n.
$$ 

For each $\nu \in {\mathbb N}$ and $m \in {\mathbb Z}_+$ define the space
$$
GS_m(M_{\nu}) = \{f \in C^m({\mathbb R}^n): q_{m, \nu}(f) = 
\sup_{x \in {\mathbb R}^n, \atop \vert \alpha \vert \le m} 
\frac {\vert (D^{\alpha}f)(x) \vert} {e^{-M_{\nu}(x)}} < \infty \}.
$$
Then, let
$GS(M_{\nu}) = \bigcap \limits_{m \in {\mathbb Z_+}}GS_m(M_{\nu})$.  
Endow $GS(M_{\nu})$ with the topology defined by the family of norms $q_{m, \nu}$ ($m \in {\mathbb Z}_+$). 
Let 
$GS({\mathcal M})= 
\bigcup \limits_{\nu \in {\mathbb N}}GS(M_{\nu})$. Supply $GS({\mathcal M})$ with an inductive limit topology of spaces $GS(M_{\nu})$. 
Let ${\varPhi^*} = \{{\varphi_{\nu}^*}\}_{\nu=1}^{\infty}$. 

In section 6 the following theorem is proved.

\begin{theorem}
Let functions of the family $\varPhi$ be convex.
Then $G({\mathcal H}) = GS({\varPhi^*})$.
\end{theorem}

The proof of Theorem 4 shows that the following more general theorem holds. 

\begin{theorem}
Let convex functions $\theta_{\nu}\colon {\mathbb R}^n \to {\mathbb R}$ 
($\nu \in {\mathbb N}$) be such that for each $\nu \in {\mathbb N}$:

1) $\theta_{\nu}(x) = \theta_{\nu}(\vert x_1 \vert, \ldots , \vert x_n \vert), \ x = (x_1, \ldots , x_n)\in {\mathbb R}^n$;

2) the restriction of $\theta_{\nu}$ to $[0, \infty)^n$ is nondecreasing in each variable;

3) there exists a constant $a_{\nu} > 0$ such that 
$\vert \varphi_{\nu} (x) - \theta_{\nu}(x)\vert \le a_{\nu}$.

Then $G({\mathcal H}) = GS({\varPhi^*})$.
\end{theorem}

In section 7 it is shown that the space $G({\mathcal H})$ represents a generalization of the Gelfand-Shilov space of type $W_M$.

\section{Auxiliary results}

In the proofs of the main results we need the following assertions. 

\begin{lemma} 
For any $m \in {\mathbb N}$, $a > 0$ there exists a constant $b_{m, a} > 0$ such that
$$
h_{m + 1}^*(x) \ge h_m^*(x) + \langle x, \tilde a \rangle - b_{m, a}, \ x \in {\mathbb R}_+^n.
$$
\end{lemma}

{\bf Proof}. Let  $a > 0$ and $x \in {\mathbb R}_+^n$. Then 
$$
h_{m + 1}^*(x) = \displaystyle \sup \limits_{y \in {\mathbb R}^n}(\langle x, y \rangle - h_{m + 1}(y)) = 
\displaystyle \sup \limits_{y \in {\mathbb R}_+^n}(\langle x, y \rangle - h_{m + 1}(y)) \ge 
$$
$$
\ge
\displaystyle \sup \limits_{y \ge \tilde a}(\langle x, y \rangle - h_{m + 1}(y)) = 
\displaystyle \sup \limits_{y \in {\mathbb R}_+^n} (\langle x, y + \tilde a \rangle - h_{m + 1}(y +\tilde a)) =
$$
$$
= \langle x, \tilde a \rangle + \displaystyle \sup \limits_{y \in {\mathbb R}_+^n} 
(\langle x, y \rangle - h_m(y) + h_m(y) - h_{m + 1}(y +\tilde a)).
$$
Further, using the condition 6) on ${\mathcal H}$ we have that
$$
h_{m + 1}^*(x) \ge \langle x, \tilde a \rangle + \displaystyle \sup \limits_{y \in {\mathbb R}_+^n} 
(\langle x, y \rangle - h_m(y))  - h_m(\tilde a) -  l_m = 
$$
$$
=\langle x, \tilde a \rangle + 
\displaystyle \sup \limits_{y \in {\mathbb R}^n} 
(\langle x, y \rangle - h_m(y))  - h_m(\tilde a) -  l_m.
$$
Hence, letting $b_{m, a} = h_m(\tilde a) + l_m$ we obtain our assertion. \ $\square$

\begin{corollary}
For any $m \in {\mathbb N}$, $a > 0$ 
$$
\varphi_m(x) + a \ln (1 + \Vert x \Vert) \le \varphi_{m+1}(x) + b_{m, a} > 0, \ x \in {\mathbb R}^n,
$$
where $b_{m, a}$ is a constant from Lemma 1.
\end{corollary}

\begin{lemma} 
For each $m \in {\mathbb N}$
$$
h_m^* (r + w) \le h_{m + 1}^*(r) + \gamma_m, \ r \in {\mathbb R}_+^n, w \in [0, \ln 2]^n.
$$
\end{lemma}

{\bf Proof}.  Let $m \in {\mathbb N}$, $r \in {\mathbb R}_+^n, w \in [0, \ln 2]^n$. Then
$$
h_m^* (r + w) = \sup_{t \in {\mathbb R}^n} (\langle r + w, t \rangle - h_m(t)) = 
\sup_{t \in {\mathbb R}_+^n}(\langle r + w, t \rangle - h_m(t)) = 
$$
$$
= \sup_{t \in {\mathbb R}_+^n}(\langle r, t \rangle + \langle w, t \rangle - h_m(t) + 
h_{m + 1}(t) - h_{m + 1}(t)).
$$
Using the condition 5) on ${\mathcal H}$ we have that 
$$
h_m^* (r + w) \le  \sup_{t \in {\mathbb R}_+^n}(\langle r, t \rangle - h_{m + 1}(t)) + \gamma_m =
$$
$$
= \sup_{t \in {\mathbb R}^n} (\langle r, t \rangle - h_{m + 1}(t)) + \gamma_m = h_{m + 1}^*(r) + \gamma_m. \ \square
$$

\begin{lemma} 
For each $m \in {\mathbb N}$ and all  $x = (x_1, \ldots , x_n) \in {\mathbb R}_+^n$ 
$$
\inf \limits_{t = (t_1, \ldots , t_n) \in (0, \infty)^n} (- x_1 \ln t_1 - \cdots - x_n \ln t_n + 
\varphi_m(t)) \le
-h_{m + 1} (x) + \gamma_m.
$$
\end{lemma}

{\bf Proof}. Using Lemma 2 for any  $x = (x_1, \ldots , x_n) \in {\mathbb R}_+^n$ we have that
$$
\inf \limits_{t = (t_1, \ldots , t_n) \in (0, \infty)^n} (- x_1 \ln t_1 - \cdots - x_n \ln t_n + 
\varphi_m(t))  =
$$
$$
= \inf \limits_{t = (t_1, \ldots , t_n) \in (0, \infty)^n} (- x_1 \ln t_1 - \cdots - x_n \ln t_n + 
h_m^*(\ln (1 + t_1), \ldots , \ln (1 + t_n)) \le 
$$
$$
\le \inf \limits_{t = (t_1, \ldots , t_n) \in [1, \infty)^n} (- x_1 \ln t_1 - \cdots - x_n \ln t_n + 
h_m^*(\ln (1 + t_1), \ldots , \ln (1 + t_n)) =
$$
$$
= - \sup_{u \in {\mathbb R}_+^n}(\langle x, u \rangle - h_m^*(\ln (1 + e^{u_1}), \ldots , \ln (1 + e^{u_n})) \le
$$
$$
\le - \sup_{u \in{\mathbb R}_+^n}(\langle x, u \rangle - h_m^*(\ln (2 e^{u_1}), \ldots , \ln (2 e^{u_n})) =
$$
$$
= - \sup_{u \in{\mathbb R}_+^n}(\langle x, u \rangle - h_m^*(\ln 2 + u_1), \ldots , \ln 2 + u_n)) \le
$$
$$
\le - \sup_{u \in{\mathbb R}_+^n}(\langle x, u \rangle - h_{m + 1}^*(u)) + \gamma_m = 
$$
$$
= - \sup_{u \in{\mathbb R}^n}(\langle x, u \rangle - h_{m + 1}^*(u)) + \gamma_m= 
-h_{m + 1} (x) + \gamma_m. \ \square
$$

\begin{corollary}
For each $m \in {\mathbb N}$ 
$$
h_{m + 1} (x) \le (\varphi_m[e])^*(x) + \gamma_m, \ x \in {\mathbb R}_+^n.
$$
\end{corollary}

\begin{lemma} 
For each $m \in {\mathbb N}$  
$$
\sup \limits_{t \in {\mathbb R}_+^n} (\langle x, t \rangle - \varphi_m [e](t)) \ge 
\sup \limits_{t \in {\mathbb R}^n} (\langle x, t \rangle - \varphi_{m + 1}  [e](t)) - \gamma_m, \ x \in {\mathbb R}_+^n.
$$
\end{lemma}

{\bf Proof}. Let $m \in {\mathbb N}$ be arbitrary. 
Consider $t = (t_1, \ldots , t_n) \in {\mathbb R}^n \setminus {\mathbb R}_+^n$. 
For simplicity assume that first $k$ coordintes of $t$ are negative and if $k < n$ then numbers 
$t_{k+1}, \ldots , t_n$ are nonnegative. 
Let $\bar t = (0, \ldots  0, t_{k+1}, \ldots , t_n)$.   
Then 
$$
\langle x, t \rangle - \varphi_{m + 1} [e](t) \le 
\langle x, \bar t \rangle - \varphi_{m + 1} [e](t).
$$
Using Lemma 2 we have that
$$
h_m^*(\ln2, \ldots , \ln 2, \ln (1 + e^{t_{k+1}}), \ldots , \ln (1 + e^{t_n})) - 
\varphi_{m + 1} [e](t) \le \gamma_m.
$$
Hence,
$$
\langle x, t \rangle - \varphi_{m + 1} [e](t) \le 
$$
$$
\le 
\langle x, \bar t \rangle - h_m^*(\ln2, \ldots , \ln 2, \ln (1 + e^{t_{k+1}}), \ldots , \ln (1 + e^{t_n})) + \gamma_m = 
$$
$$
=\langle t, \bar x \rangle  - \varphi_m [e](\bar x) + \gamma_m.
$$
This means that 
$$
\sup \limits_{t \in {\mathbb R}^n \setminus {\mathbb R}_+^n} (\langle x, t \rangle - \varphi_{m + 1}  [e](t)) \le 
\sup \limits_{t \in {\mathbb R}_+^n} (\langle x, t \rangle - \varphi_m [e](t)) + \gamma_m.
$$ 
From this the assertion of Lemma follows. \ $\square$

\begin{lemma} 
For each $m \in {\mathbb N}$  
$$
h_m(x) \ge (\varphi_{m + 1}[e])^*(x) - \gamma_m, \ x \in {\mathbb R}_+^n.
$$
\end{lemma}

{\bf Proof}. Obviously, for each $m \in {\mathbb N}$
$$
\varphi_m [e](t) \ge h_m^*(t), \ t \in {\mathbb R}_+^n.
$$
Hence, with help of Lemma 4 for $x \in {\mathbb R}_+^n$ we have that 
$$
h_m(x) =  \sup \limits_{t \in {\mathbb R}^n} (\langle x, t \rangle - h_m^*(t)) = 
\sup \limits_{t \in {\mathbb R}_+^n} (\langle x, t \rangle - h_m^*(t)) \ge 
$$
$$
\ge
\sup \limits_{t \in {\mathbb R}_+^n} (\langle x, t \rangle - \varphi_m [e](t)) \ge 
\sup \limits_{x \in {\mathbb R}^n} (\langle t, x \rangle - \varphi_{m + 1}  [e](t)) - \gamma_m \ge 
$$
$$
\ge
(\varphi_{m + 1}[e])^*(x) - \gamma_m. \ \square
$$

Note that using Proposition 6 from \cite {UMJ} we have as in Proposition 4 from \cite {CO} the following 

\begin{theorem A} 
Let $g \in C({\mathbb R}^n)$ be a convex functions satisfying the following  conditions:

1) $g(x) = g(\vert x_1 \vert, \ldots , \vert x_n \vert), \ x = (x_1, \ldots , x_n)\in {\mathbb R}^n$;

2) the restriction of $g$ to $[0, \infty)^n$ is nondecreasing in each variable;

3) $\displaystyle \lim_{x \to \infty} \frac {g(x)}{\Vert x \Vert}= + \infty$.

Then
$$
(u[e])^*(x) + (u^*[e])^*(x) =  \sum \limits _{1 \le j \le n: \atop x_j \ne 0} (x_j \ln x_j - x_j), \,
x = (x_1, \ldots , x_n) \in {\mathbb R}_+^n \setminus \{0\};
$$
$$
(u[e])^*(0) + (u^*[e])^*(0) = 0.
$$
\end{theorem A}

Fron Theorem A we have the following

\begin{corollary}
Let $u$ satisfies the conditions of Theorem A. Then 
$$
(u[e])^*(x) + (u^*[e])^*(x) \ge \sum \limits _{j=1}^n (x_j \ln (x_j  + 1) - x_j) - n, \ 
x = (x_1, \ldots , x_n) \in {\mathbb R}_+^n.
$$
\end{corollary}

The proof of the following Lemma is given in \cite {CO}.

\begin{lemma} 
Let $g = (g_1, \ldots , g_n)$ be a vector-function on ${\mathbb R}^n$ with convex components 
$g_j: {\mathbb R}^n \to [0, \infty)$ and a function $f: {\mathbb R}^n \to {\mathbb R}$ 
be such that $f_{|[0, \infty)^n}$ is convex and nondecreasing in each argument. Then $f \circ g$ is convex on ${\mathbb R}^n$.
\end{lemma}

\begin{lemma}  
For each $m \in {\mathbb N}$
$$
\varphi_{m+1}^*(\xi) \le \varphi_m^*\left(\frac {\xi}{2}\right) + \gamma_m, \ \xi \in {\mathbb R}^n.
$$
\end{lemma}

{\bf Proof}. Let $y = (y_1, \ldots , y_n) \in {\mathbb R}^n$. By Lemma 2 we have that
$$
h_{m + 1}^*(\ln (1 + \vert y_1 \vert), \ldots , \ln (1 + \vert y_n \vert)) + \gamma_m \ge 
$$
$$
\ge h_m^*(\ln (1 + \vert y_1 \vert) + \ln 2, \ldots , \ln (1 + \vert y_n \vert) + \ln 2) = 
$$
$$
=h_m^*(\ln (2 + 2\vert y_1 \vert), \ldots , \ln (2 + 2\vert y_n \vert)) \ge 
h_m^*(\ln (1 + 2\vert y_1 \vert), \ldots , \ln (1 + 2\vert y_n \vert)).
$$
Hence, for any $\xi \in {\mathbb R}^n$
$$
\varphi_{m+1}^*(\xi) = \sup \limits_{y \in {\mathbb R}^n} 
(\langle \xi, y \rangle - h_{m + 1}^*(\ln (1 + \vert y_1 \vert), \ldots , \ln (1 + \vert y_n \vert)) \le
$$
$$
\le \sup \limits_{y \in {\mathbb R}^n}  
(\langle \xi, y \rangle - 
h_m^*(\ln (1 + \vert 2 y_1 \vert), \ldots , \ln (1 + \vert 2 y_n \vert)) + \gamma_m = \varphi_m^*\left(\frac {\xi}{2}\right) + \gamma_m. \ \square
$$

\section{Proof of Theorem 1}
 
First show that the linear mapping ${\mathcal A}$ is acting from ${\mathcal E}({\mathcal H})$ into $E(\varPhi)$. 
Let $f \in {\mathcal E}({\mathcal H})$. 
Then $f \in {\mathcal E}(h_{\nu})$ for some $\nu \in {\mathbb N}$.
Hence, for each $m \in {\mathbb Z}_+$ we have that 
$$
(1+ \Vert x \Vert)^m \vert (D^{\alpha}f)(x) \vert \le  \rho_{m, \nu} (f)
\alpha! e^{-h_{\nu}(\alpha)}, \ x \in {\mathbb R}^n, \alpha \in {\mathbb Z}_+^n.
$$
We will estimate a growth  of $F_f$ using this inequality  and the Taylor series expansion 
of $F_f(z)$ ($z = x+iy, x, y \in {\mathbb R}^n$) with respect to a point $x$: 
$$
F_f(z) = \displaystyle \sum_{\vert \alpha \vert \ge 0} \frac {(D^{\alpha}f)(x)}{\alpha!} (iy)^{\alpha}.
$$
For arbitrary $m \in {\mathbb Z}_+$ we have that
$$
(1 + \Vert z \Vert)^m \vert F_f(z) \vert \le \sum_{\vert \alpha \vert \ge 0} \frac {1}{\alpha!}
(1 + \Vert x \Vert)^m (1 + \Vert y \Vert)^m
\prod \limits_{j=1}^n(1 + \vert y_j \vert)^{\alpha_j} \vert (D^{\alpha}f)(x)\vert  \le 
$$
$$
\le \rho_{m, \nu} (f) (1 + \Vert y \Vert)^m
\sum_{\vert \alpha \vert \ge 0}
e^{-h_{\nu}(\alpha)}\prod \limits_{j=1}^n(1 + \vert y_j \vert)^{\alpha_j}  \le 
$$
$$
\le
\rho_{m, \nu} (f) 
(1 + \Vert y \Vert)^m
\sum_{\vert \alpha \vert \ge 0} 
\frac {\prod \limits_{j=1}^n(1 + \vert y_j \vert)^{\alpha_j}}
{e^{h_{\nu+1}(\alpha)}}e^{h_{\nu+1}(\alpha)-h_{\nu}(\alpha)}.
$$ 
Denoting the sum of the series 
$\displaystyle\sum_{\vert \alpha \vert \ge 0} e^{h_{\nu+1}(\alpha)-h_{\nu}(\alpha)}$
by $B_{\nu}$ we have that
$$
(1 + \Vert z \Vert)^m \vert F_f(z) \vert \le B_{\nu} \rho_{m, \nu} (f) (1 + \Vert y \Vert)^m 
\sup \limits_{\vert \alpha \vert \ge 0}
\frac {\prod \limits_{j=1}^n(1 + \vert y_j \vert)^{\alpha_j}}{e^{h_{\nu+1}(\alpha)}} \le 
$$
$$
\le B_{\nu} \rho_{m, \nu} (f) 
(1 + \Vert y \Vert)^m 
e^{\sup \limits_{t = (t_1, \ldots , t_n) \in {\mathbb R}_+^n}
(t_1 \ln (1 + \vert y_1 \vert) + \cdots + t_n \ln (1 + \vert y_n \vert) -  h_{\nu+1}(t))} =
$$
$$
= B_{\nu} \rho_{m, \nu} (f) 
(1 + \Vert y \Vert)^m 
e^{\sup \limits_{t = (t_1, \ldots , t_n) \in {\mathbb R}^n}
(t_1 \ln (1 + \vert y_1 \vert) + \cdots + t_n \ln (1 + \vert y_n \vert) -  h_{\nu+1}(t))}.
$$ 
Hence,
$$
(1 + \Vert z \Vert)^m \vert F_f(z) \vert \le B_{\nu} \rho_{m, \nu} (f) (1 + \Vert y \Vert)^m 
e^{h_{\nu + 1}^*(\ln (1 + \vert y_1 \vert), \ldots , \ln (1 + \vert y_1 \vert))}.
$$
From this inequality with help of Corollary 1 we get that 
$$
(1 + \Vert z \Vert)^m \vert F_f(z) \vert \le K_{\nu, m}
\rho_{m, \nu} (f)  e^{\varphi_{\nu+2} (\vert Im z_1 \vert, \ldots , \vert Im z_n \vert)}, \ z \in {\mathbb C}^n,
$$
where $K_{\nu, m}$ is some positive constant.
Thus, for each $m \in {\mathbb Z}_+$ we have that 
\begin{equation}
p_{\nu + 2, m} (F_f) \le K_{\nu, m} \rho_{m, \nu} (f), \  f \in {\mathcal E}(h_{\nu}).
\end{equation}
Hence, 
$F_f \in E(\varphi_{\nu + 2})$. Thus, $F_f \in E(\varPhi)$. 

Note that by the inequality (2) ${\mathcal A}$ is continuous. 

Obviously, ${\mathcal A}$ is an injective linear mapping from ${\mathcal E}({\mathcal H})$ into $E(\varPhi)$. 

Show that ${\mathcal A}$ is surjective. Let $F \in E(\varPhi)$. 
Then $F \in E(\varphi_{\nu})$ for some $\nu \in {\mathbb N}$. 
Let $m \in {\mathbb Z_+}$, $\alpha = (\alpha_1, \ldots , \alpha_n) \in {\mathbb Z}_+^n$ and 
$x = (x_1, \ldots , x_n) \in {\mathbb R}^n$ be arbitrary. 
For $j =1, \ldots , n$ let $R_j$ be an arbitrary positive number. For $R = (R_1, \ldots , R_n)$ let 
$L_R(x)= \{\zeta = (\zeta_1, \ldots , \zeta_n) \in {\mathbb C}^n: \vert \zeta_j - x_j \vert = R_j, j=1, \ldots , n \}$.
Using Cauchy integral formula we have that
$$
(1+ \Vert x \Vert)^m (D^{\alpha}F)(x) =
\frac {\alpha! }{(2\pi i)^n} 
\displaystyle 
\int_{L_R(x)}
\frac 
{F(\zeta) (1+ \Vert x \Vert)^m \ d \zeta}
{(\zeta_1 - x_1)^{\alpha_1 +1} \cdots (\zeta_n - x_n)^{\alpha_n +1}} .
$$
From this we get that
$$
(1+ \Vert x \Vert)^m \vert (D^{\alpha}F)(x) \vert 
\le \frac {\alpha! }{(2\pi)^n} 
\displaystyle 
\int \limits_{L_R(x)}
\frac 
{(1+ \Vert x  - \zeta\Vert)^m (1+ \Vert \zeta \Vert)^m   \vert F(\zeta) \vert \, \vert d \zeta \vert}
{\vert \zeta_1 - x_1\vert^{\alpha_1 +1} \cdots \vert \zeta_n - x_n\vert^{\alpha_n +1}} \le 
$$
$$
\le \frac 
{\alpha! p_{\nu, m}(F) (1 + \Vert R \Vert)^m  e^{\varphi_{\nu}(R)}}{R^{\alpha}} .
$$
Using Corollary 1 we obtain that
$$
(1+ \Vert x \Vert)^m \vert (D^{\alpha}F)(x) \vert \le 
e^{b_{\nu, m}} \alpha! p_{\nu, m}(F)
\frac 
{e^{\varphi_{\nu + 1}(R)}}{R^{\alpha}}.
$$
Hence, 
$$
(1+ \Vert x \Vert)^m \vert (D^{\alpha}F)(x) \vert \le 
e^{b_{\nu, m}} \alpha! p_{\nu, m}(F)
\inf_{R \in (0, \infty)^n}
\frac 
{e^{\varphi_{\nu + 1}(R)}}{R^{\alpha}} =
$$
$$
= e^{b_{\nu, m}} \alpha! p_{\nu, m}(F)
\exp({-\sup_{r \in {\mathbb R}^n}
(\langle \alpha, r \rangle  - \varphi_{\nu + 1}[e](r))}) \le 
$$
$$
\le
e^{b_{\nu, m}} \alpha! p_{\nu, m}(F)
\exp({-\sup_{r \in {\mathbb R}_+^n}
(\langle \alpha, r \rangle  - \varphi_{\nu + 1}[e](r))})= 
$$
$$
= e^{b_{\nu, m}} \alpha! p_{\nu, m}(F) \exp({-\sup_{r \in {\mathbb R}_+^n}(\langle \alpha, r \rangle  - 
h_{\nu + 1}^*(\ln (1 +  e^{r_1}), \ldots , \ln (1 +  e^{r_n}))} \le 
$$
$$
\le e^{b_{\nu, m}} \alpha! p_{\nu, m}(F) \exp({-\sup_{r \in {\mathbb R}_+^n}(\langle \alpha, r \rangle  - 
h_{\nu + 1}^*(\ln (2 e^{r_1}), \ldots , \ln (2 e^{r_n}))} =
$$
$$
= e^{b_{\nu, m}} \alpha! p_{\nu, m}(F) \exp({-\sup_{r \in {\mathbb R}_+^n}(\langle \alpha, r \rangle  - 
h_{\nu + 1}^*(\ln 2 + r_1, \ldots , \ln 2 + r_n)))}.
$$ 
Now using Lemma 2 we have that for each $m \in {\mathbb Z}_+$  
$$
(1+ \Vert x \Vert)^m \vert (D^{\alpha}F)(x) \vert \le 
e^{b_{\nu, m}} \alpha! p_{\nu, m}(F) 
\exp(-\sup_{r \in {\mathbb R}_+^n}(\langle \alpha, r \rangle - h_{\nu + 2}^*(r)) + \gamma_{\nu}) = 
$$
$$
= e^{b_{\nu, m} + \gamma_{\nu}} \alpha! p_{\nu, m}(F) 
e^{-\sup \limits_{r \in {\mathbb R}^n}(\langle \alpha, r \rangle - h_{\nu + 2}^*(r))} = 
e^{b_{\nu, m} + \gamma_{\nu}} \alpha! p_{\nu, m}(F) e^{h_{\nu + 2}(\alpha)}.
$$
From this it follows that 
\begin{equation}
\rho_{m, \nu + 2}(F_{|{\mathbb R}^n}) \le e^{b_{\nu, m} + \gamma_{\nu}} p_{\nu, m}(F),  \ F \in E(\varphi_{\nu}). 
\end{equation}
Therefore, 
$F_{|{\mathbb R}^n} \in {\mathcal E}(h_{\nu+2})$.
Hence, $F_{|{\mathbb R}^n} \in {\mathcal E}({\mathcal H})$. Obviously, ${\mathcal A}(F_{|{\mathbb R}^n}) = F$. 
Note that the inequality (3)
ensures the continuity of ${\mathcal A}^{-1}$. Thus, the
mapping ${\mathcal A}$ establishes an isomorphism between the spaces ${\mathcal E}({\mathcal H})$ and 
$E(\varPhi)$. $\square$

\section{Proof of Proposition}

Take $g \in G({\mathcal H})$. 
Then $g \in G(h_{\nu})$ for some $\nu \in {\mathbb N}$. For arbitrary  $m \in {\mathbb Z}_+$ and 
for all $\alpha \in {\mathbb Z_+^n}$ with $\vert \alpha \vert \le m$, $\beta \in {\mathbb Z_+^n}$, 
$x \in {\mathbb R}^n$
we have that
$$
\vert x^{\beta}(D^{\alpha}g)(x) \vert \le 
\Vert g \Vert_{m, h_{\nu}} \beta! e^{-h_{\nu}(\beta)}.
$$
Using this inequality and the equality
$$
\vert (D^{\alpha}g)(x) \vert\prod \limits_{k=1}^n (1 + \vert x_k \vert)^{\beta_k}  = \vert (D^{\alpha}g)(x) \vert
\prod \limits_{k=1}^n \limits \sum \limits_{j_k=0}^{\beta_k} \binom {\beta_k}{j_k} \vert x_k \vert^{j_k} = 
$$
$$
=
\sum \limits_{j \in {\mathbb Z^n_+}: \atop (0, \ldots , 0) \le j \le \beta} \binom {\beta} {j} \vert x^j (D^{\alpha}g)(x) \vert 
$$
we get that 
$$
\vert (D^{\alpha}g)(x) \vert\prod \limits_{k=1}^n (1 + \vert x_k \vert)^{\beta_k} 
\le \Vert g \Vert_{m, h_{\nu}}
\sum \limits_{j \in {\mathbb Z^n_+}: \atop (0, \ldots , 0) \le j \le \beta} \binom {\beta}{j} 
j! e^{-h_{\nu}(j)}.
$$
From this (using the condition 6) on ${\mathcal H}$) we have that 
$$
\vert (D^{\alpha}g)(x) \vert\prod \limits_{k=1}^n (1 + \vert x_k \vert)^{\beta_k} \le 
\Vert g \Vert_{m, h_{\nu}} e^{-h_{\nu + 1}(\beta) + l_{\nu}}
\sum \limits_{j \in {\mathbb Z^n_+}: \atop (0, \ldots , 0) \le j \le \beta} \binom {\beta}{j} 
j!  e^{h_{\nu}(\beta - j)}.
$$
Therefore,
$$
\vert (D^{\alpha}g)(x) \vert\prod \limits_{k=1}^n (1 + \vert x_k \vert)^{\beta_k} \le 
\Vert g \Vert_{m, h_{\nu}} \beta! e^{-h_{\nu + 1}(\beta) + l_{\nu}}
\sum \limits_{j \in {\mathbb Z^n_+}: \atop (0, \ldots , 0) \le j \le \beta} 
\frac {e^{h_{\nu}(\beta - j)}}{(\beta - j)!} .
$$ 
In view of the condition 4) on ${\mathcal H}$ the series 
$\sum \limits_{j \in {\mathbb Z^n_+}} 
\frac {e^{h_{\nu}(j)}}{j!}$ is converging. 
From this and the inequality above it follows that 
for all $x = (x_1, \ldots , x_n) \in {\mathbb R}^n$, $\alpha \in {\mathbb Z_+^n}$ with $\vert \alpha \vert \le m$, $\beta = (\beta_1, \ldots , \beta_n) \in {\mathbb Z_+^n}$
\begin{equation}
\vert (D^{\alpha}g)(x) \vert \prod \limits_{k=1}^n (1 + \vert x_k \vert)^{\beta_k} \le c_1
\Vert g \Vert_{m, h_{\nu}} \beta! e^{-h_{\nu + 1}(\beta)} ,
\end{equation} 
where $c_1 = e^{l_{\nu}}\sum \limits_{j \in {\mathbb Z^n_+}} 
\frac {e^{h_{\nu}(j)}}{j!} $.

Now let
$$
f(\xi) = \frac {1}{(\sqrt {2 \pi})^n} 
\int_{{\mathbb R}^n} g(x) e^{i \langle x, \xi \rangle} \ dx, \ 
\xi \in {\mathbb R}^n.
$$
For all $\alpha = (\alpha_1, \ldots , \alpha_n), \beta = (\beta_1, \ldots , \beta_n) \in {\mathbb Z_+^n}, \xi \in {\mathbb R}^n$ we have that
$$
(i\xi)^{\beta}(D^{\alpha}f)(\xi) = \frac {(-1)^{\vert \beta \vert} }{(\sqrt {2 \pi})^n} 
\int_{{\mathbb R}^n} D^{\beta}(g(x) (ix)^{\alpha}) e^{i \langle x, \xi \rangle} \ dx .
$$
For $s=1, \ldots , n$ put $\gamma_s =\min (\alpha_s, \beta_s)$ and take $\gamma = (\gamma_1, \ldots , \gamma_n)$.
Then 
$$
(i\xi)^{\beta}(D^{\alpha}f)(\xi) = \frac {(-1)^{\vert \beta \vert} }{(\sqrt {2 \pi})^n} 
\int_{{\mathbb R}^n}
\displaystyle \sum \limits_{j \in {\mathbb Z_+^n}: j \le \gamma} \binom {\beta}{j}
 (D^{\beta - j} g)(x) (D^j (ix)^{\alpha}) 
e^{i \langle x, \xi \rangle} \ dx .
$$
From this we have that
$$
\vert \xi^{\beta}(D^{\alpha}f)(\xi) \vert \le \frac {1}{(\sqrt {2 \pi})^n}
\displaystyle \sum \limits_{j \in {\mathbb Z_+^n}: j \le \gamma}
\binom {\beta}{j}
\int_{{\mathbb R}^n} \vert(D^{\beta - j} g)(x)\vert 
\frac {\alpha!}{(\alpha - j)!} 
\vert x^{\alpha - j} \vert  \ dx \le
$$ 
$$
\le \frac {1}{(\sqrt {2 \pi})^n}
\displaystyle \sum \limits_{j \in {\mathbb Z_+^n}: j \le \gamma}
\frac {\binom {\beta}{j} \alpha!}{(\alpha - j)!}
\int \limits_{{\mathbb R}^n} \vert(D^{\beta - j} g)(x)\vert 
\prod \limits_{k=1}^n (1 + \vert x_k \vert)^{\alpha_k  - j_k + 2} 
\frac {dx}{\prod \limits_{k=1}^n (1 + \vert x_k \vert)^2} .
$$
Using (4) and denoting the element $(2, \ldots , 2) \in {\mathbb R}^n$ by $\omega$ we have that 
$$
\vert \xi^{\beta}(D^{\alpha}f)(\xi) \vert \le 
c_2 \Vert g \Vert_{\vert \beta \vert, h_{\nu}}
\displaystyle \sum \limits_{j \in {\mathbb Z_+^n}: j \le \gamma}
\frac {\binom {\beta}{j} \alpha!}{(\alpha - j)!} (\alpha - j + \omega)! e^{-h_{\nu + 1}(\alpha - j + \omega)}, 
$$
where
$c_2 = (\sqrt \frac {\pi}{2})^n c_1$.
Note that in view of conditions 2) and 5) on ${\mathcal H}$ we have that 
$m_{\nu} = \sup \limits_{x \in {\mathbb R}_+^n} (h_{\nu + 2} (x) - h_{\nu + 1} (x+w)) < \infty$. 
So 
$$
\vert \xi^{\beta}(D^{\alpha}f)(\xi) \vert \le 
c_2 e^{m_{\nu}} \Vert g \Vert_{\vert \beta \vert, h_{\nu}} (\alpha + w)!
\displaystyle \sum \limits_{j \in {\mathbb Z_+^n}: j \le \gamma}
\binom {\beta}{j} e^{-h_{\nu + 2}(\alpha - j)}. 
$$
Using the condition 6) on ${\mathcal H}$ we get that
$$
\vert \xi^{\beta}(D^{\alpha}f)(\xi) \vert \le 
c_2 e^{m_{\nu} + l_{\nu + 2}}\Vert g \Vert_{\vert \beta \vert, h_{\nu}} (\alpha + w)! e^{-h_{\nu + 3}(\alpha)}
\displaystyle \sum \limits_{j \in {\mathbb Z_+^n}: j \le \gamma}
\binom {\beta}{j} e^{h_{\nu + 2}(j)}.
$$
From this we have that 
$$
\vert \xi^{\beta}(D^{\alpha}f)(\xi) \vert \le 
c_2 e^{m_{\nu} + l_{\nu + 2}}\Vert g \Vert_{\vert \beta \vert, h_{\nu}} (\alpha + w)! e^{-h_{\nu + 3}(\alpha)} \beta! 
\displaystyle \sum \limits_{j \in {\mathbb Z_+^n}: j \le \gamma} \frac {e^{h_{\nu + 2}(j)}}{j!}.
$$
From this it follows that 
$$
\vert \xi^{\beta}(D^{\alpha}f)(\xi) \vert \le 
c_3 \Vert g \Vert_{\vert \beta \vert, h_{\nu}} \alpha! \beta! e^{2\vert \alpha \vert} e^{-h_{\nu + 3}(\alpha)},
$$
where
$c_3 = c_2 e^{m_{\nu} + l_{\nu + 2} + n} \displaystyle \sum \limits_{j \in {\mathbb Z_+^n}} \frac {e^{h_{\nu + 2}(j)}}{j!}$.
Using the condition 5) on ${\mathcal H}$ we get that
$$
\vert \xi^{\beta}(D^{\alpha}f)(\xi) \vert \le 
c_4 \Vert g \Vert_{\vert \beta \vert, h_{\nu}} \alpha! \beta! e^{-h_{\nu + 6}(\alpha)},
$$
where
$c_4 = c_3 e^{\gamma_{\nu + 3} + \gamma_{\nu + 4} + \gamma_{\nu + 5}}$.
Thus, for any $k \in {\mathbb Z}_+$ we can find a positive constant $c_5$ such that 
$$
(1 + \Vert \xi \Vert)^k(D^{\alpha}f)(\xi) \vert \le 
c_5 \Vert g \Vert_{k, h_{\nu}} \alpha! e^{-h_{\nu + 6}(\alpha)}, \  \alpha \in {\mathbb Z}_+^n.
$$
Hence, 
$$
\rho_{k, \nu + 6}(f) \le 
c_5 \Vert g \Vert_{k, h_{\nu}}, \ g \in G(h_{\nu}).
$$
From this it follows that the Fourier transformation ${\mathcal F}$ acts  
from $G({\mathcal H}$) into ${\mathcal E}({\mathcal H})$ and that ${\mathcal F}$ is continuous. $\square$

\section{Proof of Theorem 2}

By Theorem 1 and Proposition the linear mapping ${\mathcal A}{\mathcal F}$ acts from $G({\mathcal H}$) into $E(\varPhi)$ and is  injective  and continuous.

Show that the mapping ${\mathcal A}{\mathcal F}: G({\mathcal H}) \to E(\varPhi)$ is surjective. 
Let $F \in E(\varPhi)$. Then $F \in E(\varphi_{\nu})$ for some $\nu \in {\mathbb N}$. 
Put $f = {\mathcal A}^{-1}(F)$. Note that $f = F_|{_{{\mathbb R}^n}}$. 
By the inequality (3) $f \in {\mathcal E}(h_{\nu+2})$.
Let
$$
g(x) = \frac {1}{(\sqrt {2 \pi})^n}\int_{{\mathbb R}^n} f(\xi) 
e^{-i \langle x, \xi \rangle} \ d \xi, \ x \in {\mathbb R}^n.
$$
In other words, $g = {\mathcal F}^{-1}(f)$.
Obviously, for any $\eta \in {\mathbb R}^n$
$$
g(x) = \frac {1}{(\sqrt {2 \pi})^n}\int_{{\mathbb R}^n} F(\zeta) 
e^{-i \langle x, \zeta \rangle} \ d \xi, \  \zeta = \xi + i\eta.
$$
For any $\alpha$, $\beta = (\beta_1, \ldots , \beta_n) \in {\mathbb Z}_+^n$, $x, \eta \in {\mathbb R}^n$ we have that 
$$
x^{\beta} ({D^{\alpha} g)(x) = \frac {x^{\beta}}{(\sqrt {2 \pi})^n}\int_{{\mathbb R}^n}} F(\zeta) 
(-i \zeta)^{\alpha} 
e^{-i \langle x, \zeta \rangle} \ d \xi, \  \zeta = \xi + i\eta.
$$
From this equality we have that
$$
\vert x^{\beta} (D^{\alpha} g)(x)\vert \le \frac {1}{(\sqrt {2 \pi})^n}
\int_{{\mathbb R}^n}
\vert F(\zeta) \vert
\Vert \zeta \Vert^{\vert \alpha \vert} 
e^{\langle x, \eta \rangle} \prod \limits_{j=1}^n \vert x_j \vert^{\beta_j} \ d \xi  \le
$$
$$
\le \frac {1}{(\sqrt {2 \pi})^n}
\int_{{\mathbb R}^n} 
\vert F(\zeta) \vert
(1 + \Vert \zeta \Vert)^{n + \vert \alpha \vert + 1} 
e^{\langle x, \eta \rangle} \prod \limits_{j=1}^n \vert x_j \vert^{\beta_j} \ 
\frac {d \xi}{(1 + \Vert \xi \Vert)^{n+1}} \ .
$$
Hence, 
\begin{equation}
\vert x^{\beta} (D^{\alpha} g)(x)\vert
\le d_n p_{\nu, n + \vert \alpha \vert + 1}(F) e^{\varphi_{\nu}(\eta)}
e^{\langle x, \eta \rangle} \prod \limits_{j=1}^n \vert x_j \vert^{\beta_j} ,
\end{equation}
where $d_n = \frac {s_n}{(\sqrt {2 \pi})^n} = 
\frac 
{n}
{(\sqrt {2})^n \Gamma (\frac {n}{2} + 1)}$.
If $\vert \beta \vert = 0$ then (putting $\eta = 0$ in (5)) we have that
\begin{equation}
\vert (D^{\alpha}g)(x) \vert 
\le  d_n e^{\varphi_{\nu} (0)} p_{\nu, n + \vert \alpha \vert + 1}(F).
\end{equation} 
Now let consider the case when $\vert \beta \vert > 0$.
For $x =(x_1, \ldots , x_n) \in {\mathbb R}^n$ let $\theta(x)$ be a point in ${\mathbb R}^n$ with coordinates $\theta_j$ defined as follows: $\theta_j=\frac {x_j} {\vert x_j\vert}$ if $x_j \ne 0$ and $\theta_j=0$ if $x_j=0$ ($j=1, \ldots , n)$.
Let $t = (t_1, \ldots, t_n) \in {\mathbb R}^n$ have strictly positive coordinates.
Put
$\eta = -(\theta_1 t_1, \ldots , \theta_n t_n)$.
Then  from (5) we get that
$$
\vert x^{\beta} 
(D^{\alpha} g)(x)\vert \le d_n
 p_{\nu, n + \vert \alpha \vert + 1}(F) 
e^{\varphi_{\nu}(t)} \prod \limits_{j \in \{1, \ldots , n\}: \beta_j \ne 0} 
\frac {\vert x_j \vert^{\beta_j}}{e^{t_j \vert x_j \vert}}  \le
$$
$$
\le d_n
 p_{\nu, n + \vert \alpha \vert + 1}(F) 
e^{\varphi_{\nu}(t)} \prod \limits_{j \in \{1, \ldots , n\}: \beta_j \ne 0} 
\sup \limits_{r_j > 0}
\frac {r_j^{\beta_j}}{e^{t_j r_j}} =
$$
$$
= d_n
 p_{\nu, n + \vert \alpha \vert + 1}(F)
\exp (\varphi_{\nu}(t) + \sum \limits_{1 \le j \le n: \beta_j \ne 0} \sup \limits_{r_j > 0}(- t_j r_j + \beta_j \ln r_j)) =
$$
$$
= d_n
 p_{\nu, n + \vert \alpha \vert + 1}(F)\exp (\varphi_{\nu}(t) + \sum \limits_{1 \le j \le n: \beta_j \ne 0} (\beta_j \ln \beta_j - \beta_j) - \sum \limits_{j=1}^n \beta_j \ln t_j). 
$$
Hence,
$$
\vert x^{\beta} 
(D^{\alpha} g)(x)\vert \le d_n
p_{\nu, n + \vert \alpha \vert + 1}(F) 
e^{\sum \limits_{1 \le j \le n: \atop \beta_j \ne 0} \beta_j \ln \frac {\beta_j}{e}  + 
\inf \limits_{t = (t_1, \ldots , t_n) \in (0, \infty)^n} 
(\varphi_{\nu}(t) - \sum \limits_{j=1}^n \beta_j \ln t_j)}. 
$$
From this with help of Lemma 3 we get that 
$$
\vert x^{\beta} 
(D^{\alpha} g)(x)\vert \le d_n
 p_{\nu, n + \vert \alpha \vert + 1}(F) 
e^{\sum \limits_{1 \le j \le n: \atop \beta_j \ne 0} \beta_j \ln \frac {\beta_j}{e}  - h_{\nu + 1} (\beta) + \gamma_{\nu}}. 
$$
Hence, for $\vert \beta \vert > 0$
$$
\vert x^{\beta} 
(D^{\alpha} g)(x)\vert \le d_{n, \nu} p_{\nu, n + \vert \alpha \vert + 1}(F) \beta! e^{- h_{\nu + 1} (\beta)}, 
$$
where $d_{n, \nu}  = d_n e^{\gamma_{\nu}}$.
Now from this and (6) we have that for $\alpha, \beta \in {\mathbb Z}_+^n$, 
$x \in {\mathbb R}^n$
$$
\vert x^{\beta} 
(D^{\alpha} g)(x)\vert \le d_{n, \nu} p_{\nu, n + \vert \alpha \vert + 1}(F) \beta! e^{- h_{\nu + 1} (\beta)}.
$$
From this it follows that for each $m \in {\mathbb Z}_+$
$$
\max_{\vert \alpha \vert \le m} \sup_{x \in {\mathbb R}^n, \beta \in {\mathbb Z}^n_+}
\frac {\vert x^{\beta} 
(D^{\alpha} g)(x)\vert}{\beta! e^{- h_{\nu + 1} (\beta)}} \le 
d_{n, \nu} p_{\nu, n + m + 1}(F), \ F \in E(\varphi_{\nu}).
$$
This means that for each $m \in {\mathbb Z}_+$
$$
\Vert g \Vert_{m, h_{\nu + 1}} \le d_{n, \nu}  p_{\nu, n + m + 1}(F), \ F \in E(\varphi_{\nu}).
$$
In other words, for each $m \in {\mathbb Z}_+$
$$
\Vert {\mathcal F}^{-1} {\mathcal A}^{-1} (F) \Vert_{m, h_{\nu + 1}} \le d_{n, \nu}  p_{\nu, n + m + 1}(F), \ F \in E(\varphi_{\nu}).
$$
From this inequality it follows that the linear mapping ${\mathcal A} {\mathcal F}$ is surjective and the linear mapping
$({\mathcal A} {\mathcal F})^{-1}$ acts from $E(\varPhi)$ to $G({\mathcal H})$ and is continuous. Hence, the mapping
${\mathcal A} {\mathcal F}$ establishes an isomorphism between  $G({\mathcal H})$ and $E(\varPhi)$. $\square$

\section{A special case of ${\mathcal H}$}

{\bf Proof of Theorem 4}. Let $f \in G({\mathcal H})$. Then $f \in G(h_{\nu})$ for some  $\nu \in {\mathbb N}$. 
Let $m \in {\mathbb Z}_+$ be arbitrary. 
For all $\alpha \in {\mathbb Z_+^n}$ with $\vert \alpha \vert \le m$,  
$\beta  \in {\mathbb Z}_+^n$ and 
$x = (x_1, \ldots , x_n) \in {\mathbb R}^n$ with non-zero coordinates we have that
$$
\vert (D^{\alpha} f)(x) \vert \le 
\frac 
{\Vert f \Vert_{m, h_{\nu}}
\beta! 
e^{-h_{\nu}(\beta)}}
{\prod \limits_{j=1}^n \vert x_j \vert^{\beta_j}} \ .
$$
Take into account that 
$j! \le \frac {(j+1)^{j+1}}{e^j}$ for all $j \in {\mathbb Z_+}$. 
Then for all  
$\beta  \in {\mathbb Z}_+^n$ and 
$x = (x_1, \ldots , x_n) \in {\mathbb R}^n$
\begin{equation}
\vert (D^{\alpha} f)(x) \vert \le \Vert f \Vert_{m, h_{\nu}}  
e^{-h_{\nu}(\beta)} 
{\prod \limits_{j=1}^n \frac {(\beta_j + 1)^{\beta_j + 1}}
{(e\vert x_j \vert)^{\beta_j}}} \ .
\end{equation}
Our aim is to obtain a suitable estimate of 
$e^{-h_{\nu}(\beta)} 
{\prod \limits_{j=1}^n \frac {(\beta_j + 1)^{\beta_j + 1}}
{(e\vert x_j \vert)^{\beta_j}}}$ from above.
For $\beta \in {\mathbb Z}_+^n$ let 
$\Omega_{\beta} = \{t = (t_1, \ldots , t_n) \in {\mathbb R}^n: 
\beta_j \le t_j  < \beta_j + 1 \ (j =1, \ldots , n) \}$. 
Also, for $\lambda > 0$ let $\check \lambda:  = \max (\lambda, 1)$.
Using Lemma 2 and nondecreasity in each variable of $h$  in ${\mathbb R}_+^n$ 
for $t \in \Omega_{\beta}$ and $\mu = (\mu_1, \ldots , \mu_n) \in (0, \infty)^n $
we have that   
$$
e^{-h_{\nu}(\beta)} 
{\prod \limits_{j=1}^n \frac {(\beta_j + 1)^{\beta_j + 1}}
{\mu_j^{\beta_j}}} 
\le e^{-h_{\nu + 1}(t) + h_{\nu}(1, \ldots , 1) + l_{\nu}}
{\prod \limits_{j=1}^n \frac {\check \mu_j (t_j + 1)^{t_j + 1}}
{\mu_j^{t_j}}}.
$$
From this we get that
$$
\inf_{\beta \in {\mathbb Z}_+^n} 
e^{-h_{\nu}(\beta)} {\prod \limits_{j=1}^n \frac {(\beta_j + 1)^{\beta_j + 1}}{\mu_j^{\beta_j}}} 
\le C_1 
e^{(\inf \limits_{t \in {\mathbb R}_+^n}\sum \limits_{j=1}^n (\ln \check \mu_j  + 
(t_j + 1) \ln (t_j + 1) - t_j \ln \mu_j) -h_{\nu + 1}(t))},
$$
where
$C_1 = e^{h_{\nu}(1, \ldots , 1) + l_{\nu}}$, $t = (t_1, \ldots , t_n)$.
Now with help of Lemma 5 we get 
$$
\inf_{\beta \in {\mathbb Z}_+^n} 
e^{-h_{\nu}(\beta)} {\prod \limits_{j=1}^n \frac {(\beta_j + 1)^{\beta_j + 1}}{\mu_j^{\beta_j}}} 
\le 
$$
$$
\le C_2 
e^{\inf \limits_{t \in {\mathbb R}_+^n} (\sum \limits_{j=1}^n (\ln \check \mu_j  + 
(t_j + 1) \ln (t_j + 1) - t_j \ln \mu_j) - (\varphi_{\nu + 2} [e])^*(t))},
$$
where $C_2 = C_1 e^{\gamma_{\nu + 1}}$.
Now using this Corollary 3 we have that
$$
e^{-h_{\nu}(\beta)} 
{\prod \limits_{j=1}^n \frac {(\beta_j + 1)^{\beta_j + 1}}
{\mu_j^{\beta_j}}} 
\le C_3 
e^{\sum \limits_{j=1}^n (\ln \check \mu_j  + \ln (t_j + 1) - t_j \ln \mu_j + t_j) + 
(\varphi_{\nu + 2}^*[e])^*(t)},
$$
where $C_3 = C_2 e^n$.
Obviously, there exists a constant $C_4 > 0$ such that 
$$
e^{-h_{\nu}(\beta)} 
{\prod \limits_{j=1}^n \frac {(\beta_j + 1)^{\beta_j + 1}}
{\mu_j^{\beta_j}}}  \le 
C_4 
e^{(\varphi_{\nu + 2}^*[e])^*(t) - \sum \limits_{j=1}^n  t_j \ln \frac {\mu_j}{4} + 
\sum \limits_{j=1}^n \ln \check \mu_j}.
$$
From this it follows that
$$
\inf_{\beta \in {\mathbb Z}_+^n}
e^{-h_{\nu}(\beta)}  
{\prod \limits_{j=1}^n \frac {(\beta_j + 1)^{\beta_j + 1}}
{\mu_j^{\beta_j}}}  \le 
C_4 
e^{\inf \limits_{t \in [0, \infty)^n} ((\varphi_{\nu + 2}^*[e])^*(t) - \sum \limits_{j=1}^n  t_j \ln \frac {\mu_j}{4}) + 
\sum \limits_{j=1}^n \ln \check \mu_j}.
$$
In other words,
$$
\inf_{\beta \in {\mathbb Z}_+^n}
e^{-h_{\nu}(\beta)} 
{\prod \limits_{j=1}^n \frac {(\beta_j + 1)^{\beta_j + 1}}
{\mu_j^{\beta_j}}}  \le 
C_4
e^{-\sup \limits_{t \in [0, \infty)^n} (\sum \limits_{j=1}^n  t_j \ln \frac {\mu_j}{4} - 
(\varphi_{\nu + 2}^*[e])^*(t)) + 
\sum \limits_{j=1}^n \ln \check \mu_j}.
$$
Taking into account that the function  
$(\varphi_{\nu + 2}^*[e])^*$ takes finite values on $[0, \infty)^n$ and $(\varphi_{\nu + 2}^*[e])^*(x) = +\infty$ if 
$x \notin [0, \infty)^n$ we can rewrite the last inequality in the following form
$$
\inf_{\beta \in {\mathbb Z}_+^n}
e^{-h_{\nu}(\beta)} 
{\prod \limits_{j=1}^n \frac {(\beta_j + 1)^{\beta_j + 1}}
{\mu_j^{\beta_j}}}  \le 
C_4 
e^{-\sup \limits_{t \in {\mathbb R}^n} (\sum \limits_{j=1}^n  t_j \ln \frac {\mu_j}{4} - 
(\varphi_{\nu + 2}^*[e])^*(t)) + 
\sum \limits_{j=1}^n \ln \check \mu_j}.
$$
Note that by Lemma 6 the function $\varphi_{\nu + 2}^*[e]$ is convex on ${\mathbb R}^n$ with finite values (thus, $\varphi_{\nu + 2}^*[e]$ is continuous on ${\mathbb R}^n$ (see \cite {R}, Corollary 10.1.1)).
Taking into account that the Young-Fenchel conjugation is involutive (see \cite {R}, Theorem 12.2) we have that 
$$
\sup \limits_{t \in {\mathbb R}^n} (\sum \limits_{j=1}^n  t_j \ln \frac {\mu_j}{4} - 
(\varphi_{\nu + 2}^*[e])^*(t)) = 
\varphi_{\nu + 2}^*[e]
\left(\ln \frac {\mu_1}{4}, \ldots , \ln \frac {\mu_n}{4}\right) .
$$
Thus,
$$
\inf_{\beta \in {\mathbb Z}_+^n}
e^{-h_{\nu}(\beta)}  
{\prod \limits_{j=1}^n \frac {(\beta_j + 1)^{\beta_j + 1}}
{\mu_j^{\beta_j}}}  \le 
C_4 
e^{-\varphi_{\nu + 2}^*[e]
\left(\ln \frac {\mu_1}{4}, \ldots , \ln \frac {\mu_n}{4}\right) + 
\sum \limits_{j=1}^n \ln \check \mu_j}.
$$
In other words, 
$$
\inf_{\beta \in {\mathbb Z}_+^n}
e^{-h_{\nu}(\beta)} 
{\prod \limits_{j=1}^n \frac {(\beta_j + 1)^{\beta_j + 1}}
{\mu_j^{\beta_j}}}  \le 
C_4 
e^{-\varphi_{\nu + 2}^*(\frac {\mu}{4}) + 
\sum \limits_{j=1}^n \ln \check \mu_j}.
$$
From this inequality with help of Lemma 7 we get that 
\begin{equation}
\inf_{\beta \in {\mathbb Z}_+^n}
e^{-h_{\nu}(\beta)}
{\prod \limits_{j=1}^n \frac {(\beta_j + 1)^{\beta_j + 1}}
{\mu_j^{\beta_j}}}  \le 
C_4 
e^{-\varphi_{\nu + 4}^*(\mu) + 
\sum \limits_{j=1}^n \ln \check \mu_j},
\end{equation}
where $C_4 = C_3 e^{a_{\nu + 1} + \gamma_{\nu + 2} + \gamma_{\nu + 3}}$.
Note that by Lemma 7 and Lemma 5 in \cite {CO} for each $j \in {\mathbb N}$ 
we have that
$$
\displaystyle \lim_{\xi \to \infty} \frac 
{\varphi_j^*(\xi) - \varphi_{j+1}^*(\xi)}{\Vert \xi \Vert}= + \infty.
$$
Using this we obtain from (8) that
$$
\inf_{\beta \in {\mathbb Z}_+^n}
e^{-h_{\nu}(\beta)} 
{\prod \limits_{j=1}^n \frac {(\beta_j + 1)^{\beta_j + 1}}
{\mu_j^{\beta_j}}}  \le 
C_5 
e^{-\varphi_{\nu + 5}^*(\mu)},
$$
where
$C_5$ is some positive number depending on $\nu$.
From this and the inequality (7) we obtain at last that for all 
$x = (x_1, \ldots , x_n) \in {\mathbb R}^n$ with non-zero coordinates and for all 
$\alpha \in {\mathbb Z_+^n}$ with $\vert \alpha \vert \le m$
$$
\vert (D^{\alpha} f)(x) \vert \le  C_5 \Vert f \Vert_{m, h_{\nu}}
e^{-\varphi_{\nu + 5}^*(e \vert x_1 \vert, \ldots , e \vert x_n \vert)} \le
C_5 \Vert f \Vert_{m, h_{\nu}}
e^{-\varphi_{\nu + 5}^*(x)}.
$$
Obviously, the last inequality holds for all $x \in {\mathbb R}^n$.  
Thus,  $f \in GS(\varphi_{\nu + 5}^*)$ and
$$
q_{m, \nu + 5}(f) \le C_5 \Vert f \Vert_{m, h_{\nu}}, \ f \in G(h_{\nu}).
$$
From this it follows that the identity mapping $J$ acts from $G({\mathcal H})$ to $GS({\varPhi^*})$ and is continuous.

Show that $J$ is surjective. Let $f \in GS(\varPhi^*)$. 
Then $f \in GS(\varphi_{\nu}^*)$ for some $\nu \in {\mathbb N}$. 
Fix $m \in {\mathbb Z}_+$. 
Consider an arbitrary point $x = (x_1, \ldots , x_n) \in {\mathbb R}^n$  with non-zero coordinates.
For all $\alpha \in {\mathbb Z_+^n}$ with $\vert \alpha \vert \le m$ we have that
$$
\vert (D^{\alpha}f)(x) \vert \le q_{m, \nu}(f)
e^{-\varphi_{\nu}^*(\vert x_1 \vert, \ldots , \vert x_n \vert)}.
$$
In other words,
$$
\vert (D^{\alpha}f)(x) \vert   \le q_{m, \nu}(f)
e^{-\varphi_{\nu}^*[e](\ln \vert x_1 \vert, \ldots , \ln \vert x_n \vert)}.
$$
From this we have that
$$
\vert (D^{\alpha}f)(x) \vert  \le q_{m, \nu}(f)  
\exp (- \sup \limits_{t =(t_1, \ldots , t_n) \in {\mathbb R}_+^n} 
(\sum \limits_{j=1}^n t_j \ln \vert x_j \vert - (\varphi_{\nu}^*[e])^*(t))).
$$
Now with help of Theorem A we get that
$$
\vert (D^{\alpha}f)(x) \vert \le q_{m, \nu}(f)  
e^{- \sup \limits_{t =(t_1, \ldots , t_n) \in {\mathbb R}_+^n} 
(\sum \limits_{1 \le j \le n:  t_j \ne 0} (t_j \ln (e \vert x_j \vert)  - t_j \ln t_j)+ 
(\varphi_{\nu}[e])^*(t))}.
$$
Using Corollary 2 we get that 
$$
\vert (D^{\alpha}f)(x) \vert \le q_{m, \nu}(f)  e^{\gamma_{\nu}}
e^{- \sup \limits_{t =(t_1, \ldots , t_n) \in {\mathbb R}_+^n} 
(\sum \limits_{1 \le j \le n:  t_j \ne 0} (t_j \ln (e \vert x_j \vert)  - t_j \ln t_j)+ 
h_{\nu + 1}(t))}.
$$
Consequently, if $\beta \in {\mathbb Z}_+^n$ then
$$
\vert (D^{\alpha}f)(x) \vert \le q_{m, \nu}(f) e^{\gamma_{\nu}} e^{-h_{\nu + 1}(\beta)}
\prod \limits_{1 \le j \le n: \beta_j \ne 0} 
\frac 
{\beta_j^{\beta_j}}
{(e \vert x_j \vert)^{\beta_j}}. 
$$
From this we finally obtain that for all $x \in {\mathbb R}^n$  with non-zero coordinates
$$
\vert x^{\beta}(D^{\alpha}f)(x) \vert  \le 
q_{m, \nu}(f) e^{\gamma_{\nu}}
\beta!  e^{-h_{\nu + 1}(\beta)}, \ \beta \in {\mathbb Z}_+^n, \vert \alpha \vert \le m.
$$
Clearly, this inequality holds for all $x \in {\mathbb R}^n$. 
Thus,
$$
\Vert f \Vert_{m, h_{\nu + 1}^*} \le q_{m, \nu}(f).
$$
Since here $m \in {\mathbb Z}_+$ is arbitrary then $f \in G(h_{\nu + 1})$. Hence, $f \in G({\mathcal H})$.
Also from the last inequality it follows that the mapping $J^{-1}$ is continuous. Thus, 
the spaces $G({\mathcal H})$ and $GS({\varPhi^*})$ coincide. $\square$

\section{Appendix}

Let ${\mathcal M} = \{{M_{\nu}}\}_{\nu=1}^{\infty}$ be a family of convex functions 
$M_{\nu}: {\mathbb R}^n \to {\mathbb R}$ satisfying the  conditions $j_1)-j_4)$. 
Put $\varphi_{\nu} = M_{\nu}^*$. Then $M_{\nu} = \varphi_{\nu}^*$. 

By Lemma 6 for each $\nu \in {\mathbb N}$ the function $\psi_{\nu}$ on ${\mathbb R}^n$ defined by the formula 
$$
\psi_{\nu}(x) = \varphi_{\nu}[e](\vert x_1 \vert, \ldots , \vert x_n \vert), \ x = (x_1,   \ldots , x_n) \in {\mathbb R}^n,
$$
is convex. Obviously, ${\psi_{\nu}}_{|{\mathbb R}^n}$ is nondecreasing in each variable. 

Put $h_{\nu} = \psi_{\nu}^*, \ \nu \in {\mathbb N}$. Since 
$$
\displaystyle \lim_{x \to \infty} \frac {\psi_{\nu}(x)}{e^{\vert x_1 \vert} + \cdots + e^{\vert x_n \vert}}  = 
+ \infty,
$$
then convex function $h_{\nu}$ takes finite values and satisfies the conditions 3) and 4) (by Lemma 1 and Remark 2 in \cite {CO}). Obviously,  conditions 1) and 2) are fulfilled for $h_{\nu}$.
Further, in view of the condition $j_4)$ we have that 
$$
M_{\nu}^*(2 x) \le M_{\nu + 1}^*(x) + A_{\nu}, \ x \in {\mathbb R}^n.
$$
From this it follows that 
\begin{equation}
\varphi_{\nu}(2 x) \le \varphi_{\nu + 1}(x) + A_{\nu}, \ x \in {\mathbb R}^n.
\end{equation}
Then by Lemma 3 in \cite {CO} functions $h_{\nu}$ satisfy the conditions 5) with 
$\gamma_{\nu} = A_{\nu}$. 
At last, using convexity of functions $\varphi_{\nu}$ from (9) we get that 
$$
2 \varphi_{\nu}(x) \le \varphi_{\nu + 1}(x) + \varphi_{\nu}(0) + A_{\nu}, \ x \in {\mathbb R}^n.
$$
Now by Lemma 2 in \cite {CO} functions $h_{\nu}$ satisfy the condition 6) with 
$l_{\nu} = \varphi_{\nu}(0) + A_{\nu}$. 
Next, since functions $M_{\nu}$ are convex then in view of the condition $j_4)$ we have that 
$$
2 M_{\nu + 1}(x) \le M_{\nu}(x) + A_{\nu} + M_{\nu}(0), \ x \in {\mathbb R}^n. 
$$
Hence,
$$
M_{\nu}(x) - M_{\nu + 1}(x) \ge M_{\nu + 1}(x) - A_{\nu} - M_{\nu}(0), \ x \in {\mathbb R}^n. 
$$
Then for each $\nu \in {\mathbb N}$ there exists a constant $K_{\nu} > 0$ such that 
$$
\varphi_{\nu}(x + \xi) \le \varphi_{{\nu}+1}(x) + K_{\nu}, \ x \in [0, \infty)^n, \xi \in [0, 1]^n.
$$
This inequality and the inequality mean that functions $\varphi_{\nu}$ satisfy the conditions of Theorem 4 in \cite {CO}. So by this theorem $GS({\mathcal M}) = G({\mathcal H})$.

\vspace {0.3cm}

Institute of Mathematics with Computer Centre 

of Ufa Scientific Centre of Russian Academy of Sciences, 

Chernyshevsky str., 112, Ufa, 450077, Russia

E-mail: musin\_ildar@mail.ru


\begin{thebibliography}{99}

\bibitem {GS1} 
I.M. Gelfand, G.E. Shilov, {\it Generalized functions}, Vol. 2, Academic Press, New York,
1968.

\bibitem {GS2} 
I.M. Gelfand, G.E. Shilov, {\it Generalized functions}, Vol. 3, Academic Press, New York, 1967.

\bibitem {CO} Il'dar Kh. Musin, On a space of entire functions rapidly decreasing on ${\mathbb R}^n$ and its Fourier transform,  {\it Concrete Operators}, 2:1 (2015), 120-138. 

\bibitem {UMJ} I.Kh. Musin, On a Hilbert space of entire functions, {\it Ufa Math. Journal}, 9:3 (2017), 109-117.

\bibitem {R} R.T. Rockafellar, {\it Convex analysis},  Princeton, New Jersey, Princeton University Press, 1970.


\end{thebibliography}
\end{document}